\documentclass[a4paper,11pt]{article}

\usepackage{amsmath}
\usepackage{amssymb}
\usepackage{amsthm}
\usepackage{mathrsfs}
\usepackage{enumerate}

\newtheorem{thm}{Theorem}[section]
\newtheorem{lemma}[thm]{Lemma}

\newtheorem{defn}{Definition}[section]

\newtheorem{remark}{Remark}[section]
\newtheorem*{notation}{Notation}

\newcommand{\C}{\mathscr{C}}
\newcommand{\disth}{d_{\scriptscriptstyle{\mathrm{H}}}}
\newcommand{\eps}{\varepsilon}
\newcommand{\Pro}{\mathtt{P}}
\newcommand{\Ca}{\mathtt{C}}

\newcommand{\set}[1]{\left\{#1\right\}}
\newcommand{\R}{\mathbb R}

\newcommand{\N}{\mathbb N}

\newcommand{\proba}{\mathbb P}
\newcommand{\I}{\mathcal{I}}

\newcommand{\Com}{\mathbb C}
\newcommand{\virg}[1]{``#1"}

\begin{document}

\title{{\bf $\eps$-Distortion Complexity for Cantor Sets}}

\author{{\bf C. Bonanno}\\
{\em Dipartimento di Matematica Applicata}\\
{\em Universit\`a di Pisa}
{\em 56127 Pisa, Italy}\\  \\
{\bf J.-R. Chazottes, P. Collet}\\
{\em Centre de Physique Th\'eorique}\\
{\em  Ecole polytechnique, CNRS}\\ 
{\em 91128 Palaiseau, France}
}

\maketitle

\begin{abstract}
We define the $\eps$-distortion complexity of a set as the shortest
program, running on a universal Turing machine, which produces this set at the
precision $\eps$ in the sense of Hausdorff distance. Then, we estimate
the $\eps$-distortion complexity of various central Cantor sets on the
line generated by iterated function systems (IFS's). In particular, the
$\eps$-distortion complexity of a $C^k$ Cantor set depends, in general,
on $k$ and on its box counting dimension, contrarily to Cantor sets generated by
polynomial IFS or random affine Cantor sets.

\bigskip

\noindent {\bf keywords}: central Cantor sets, iterated function
system, random Cantor sets, scaling function, $C^k$ Cantor sets, box
counting dimension.

\end{abstract}

\section{Introduction}

Nowadays, computers are being widely used to generate images
in the analysis and simulations of real-life processes and their
mathematical models. A natural issue is to measure the complexity
of drawing a set of points on a computer, which describes a continuous
object at a given precision.
Particular examples of complex objects are fractal sets which arise in many contexts \cite{fractalbook}. 
Well-known examples of fractal sets are strange attractors of dissipative dynamical systems and Julia sets.
Another way to generate fractal sets is to use iterated function systems \cite{barnsley}.

The way we measure the complexity of a (compact) set can be
colloquially described as follows. We define the $\eps$-distortion
complexity of a set $\C$ as the minimal length of the programs
producing a finite set $\eps$-close to $\C$, in the sense of
Hausdorff distance. As in the classical notion of Kolmogorov complexity of sequences,
by programs we mean programs running on a universal Turing machine \cite{li-vit}.   
We are interested in the behavior of  the $\eps$-distortion
complexity when $\eps$ is getting small, and in eventual relations of this behavior
with other characteristics of the set ({\em e.g.}, fractal dimension).

In the present article we consider various classes of Cantor sets on the real
line generated by iterated function systems (IFS's) \cite{barnsley} and
compute bounds from above and below of their $\eps$-distortion
complexity as a function of $\eps$.
We first consider IFS's with polynomial contractions and obtain the upper bound
$const\times\log(\eps^{-1})$ for the $\eps$-distortion
complexity of the generated Cantor set, where the (finite) constant may depend on the polynomials. 
We can produce ``many'' polynomial IFS's
with a lower bound of the same order using a probabilistic construction.
It turns out  that some particular Cantor sets like the usual middle
third Cantor set are of much lower complexity. 
For analytic IFS's, we obtain the upper bound $const\times (\log(\eps^{-1}))^2$.
Next we consider random central Cantor sets produced by affine IFS's,
for which the contraction rate is chosen at random at each step of
the construction. In this case, we get the upper bound
$const\times(\log(\eps^{-1}))^2$ and the lower bound
$const\times(\log(\eps^{-1}))^{2-\delta}$,
for any $\delta>0$, for almost all such Cantor sets (where the constant in our bound depends on $\delta$
and tends to $0$ when $\delta\to 0$). Finally, we 
consider $C^k$ IFS's. Contrarily to the previous cases, the leading, asymptotic behavior
of the $\eps$-distortion complexity depends on the box counting
dimension $D$ of the generated Cantor set. Indeed, we obtain the upper
bound $const\times\eps^{-\frac{D}{k}-\delta}$, for any $\delta>0$ (where the constant in our bound depends on $\delta$
and blows up when $\delta\to 0$).
We then construct ``many'' $C^k$ (random) Cantor sets
with a lower bound $const\times\eps^{-\frac{D}{k}+\delta}$ (for any $\delta>0$),
by constructing their scaling function \cite{sullivan}.

The case of sets reduced to one point on the line was investigated in \cite{cai} where in particular the Hausdorff dimension
of the set of reals with given asymptotic complexity is computed. 
For graphs of functions, from the point of view of determining the values of a function at given precision, relations
with $\eps$-entropy are obtained in \cite{asarin}.
Another notion of complexity is to ask about the smallest execution time of the programs  generating a given set with $\eps$-precision,
in the sense of Hausdorff distance \cite{W}. This was used in \cite{braverman} to show that a class of Julia sets
was polynomial time computable.

This article is organized as follows. In Section \ref{definition}
we define the $\eps$-distortion complexity of a compact set and
state our results. Section \ref{allproofs} is devoted to the proofs.

\section{Definitions and results} \label{definition}

For a compact set $\C \subset \R^d$, we define its \emph{$\eps$-distortion
complexity} as follows.

\begin{defn} 
The $\eps$-distortion complexity of a compact set $\C \subset
\R^d$ at precision $\eps>0$ is defined by
$$
\Delta(\C,\eps)=\min\left\{ \ell(\Pro) : \disth(\Ca(\Pro),\C)<\eps
\right\}\ ,
$$
where the minimum is taken over all binary programs $\Pro \in
\set{0,1}^*$ running on a universal Turing machine $U$, which
produce a finite subset $\Ca(\Pro) \in \R^d$; $\ell(\Pro)$ is the
program length; $\disth$ denotes the Hausdorff distance.
\end{defn}

Notice that, because of the compactness of $\C$, we can
use a minimum in the above definition, which always leads to a
finite number.

For the reader's convenience, we recall that the Hausdorff distance $\disth$
between two closed subsets $F_1,F_2$ of a metric space with metric $d$ is
given by (see, {\em e.g.}, \cite{mattila,barnsley})
$$
\disth(F_1,F_2)=\max\{ \sup_{x_2\in F_2} d(x_2,F_1),\sup_{x_1\in F_1} d(x_1,F_2)\} .
$$

\begin{remark}
If $T$ is a bi-Lipschitz map, there exist two positive constants $c_1,c_2$ such that, if 
$\C$ is a compact set, we have
$$
\Delta(T(\C),c_1 \eps)\leq \Delta(\C,\eps) \leq \Delta(T(\C),c_2 \eps) .
$$
\end{remark}

We now recall the definition of Cantor sets generated by iterated
function systems \cite{barnsley}. 
For the sake of simplicity, we restrict ourselves to Cantor sets in the unit
interval $[0,1]\subset \R$, although several results can be easily
generalised to arbitrary finite dimension.

Let $A=[0,1]$ and let $\I$ be a finite set of indices with at
least two elements. An \emph{Iterated Function System} (IFS for
short) is a collection
$$
\set{\phi_i : A \to A  : i\in \I}
$$
of injective contractions on $A$ with uniform contraction rate
$\rho \in (0,1)$, and such that $\phi_i(A) \cap \phi_j(A) =
\emptyset$ for $i\not= j$.

For any infinite word $\omega \in \I^\infty$ and for any $n\in
\N$, let $\omega_1^n \in \I^n$ denote the prefix of length $n$
given by the first $n$ symbols of $\omega$, and let
\begin{equation} \label{ifs-finite}
\phi_{_{\omega_1^n}} := \phi_{\omega_n} \circ \phi_{\omega_{n-1}}
\circ \dots \circ \phi_{\omega_1}.
\end{equation}
The map $\pi:\I^\infty \to A$ defined by
$
\omega \mapsto \pi(\omega):= \bigcap_{n=0}^\infty \
\phi_{_{\omega_1^n}}(A)
$
is continuous (in product topology) and, since
$$
\hbox{diam } (\phi_{_{\omega_1^n}} (A)) \le \rho^n \hbox{ diam }(A),
$$
$\pi(\omega)$ is a point in $A$ for all $\omega \in \I^\infty$.
The set 
$$
\C:= \pi(\I^\infty) = \bigcup_{\omega \in \I^\infty}\
\bigcap_{n=0}^\infty \ \phi_{_{\omega_1^n}}(A)
$$
is a Cantor set and satisfies
\begin{equation} \label{ifs-prop}
\C = \bigcup_{i\in \I} \ \phi_i(\C).
\end{equation}

We are interested in the behaviour of $\Delta(\C,\eps)$ when
$\eps$ tends to zero. Note that this is a monotone decreasing
function.

\begin{notation}
In the sequel we write $f\asymp g$ if there are two positive
constants $C_1$ and $C_2$ such that for any $\eps>0$ small enough
$$
C_1 f(\eps)\le g(\eps)\le C_2 f(\eps).
$$
We write $f\preccurlyeq g$ if there is a positive constant $C$
such that for any $\eps>0$ small enough
$$
f(\eps)\le C g(\eps).
$$
\end{notation}

\bigskip

Our first result deals with polynomial IFS's.
\begin{thm} \label{thm:pol}
Let $\C$ be a Cantor set generated by an IFS with polynomial
functions. Then
\begin{equation} \label{pol-1}
\Delta(\C,\eps) \preccurlyeq \log(\eps^{-1}).
\end{equation}
Moreover, for any $\delta>0$, there exist (many) polynomial IFS's
such that the generated Cantor set satisfies
\begin{equation} \label{pol-2}
(1-\delta) \log(\eps^{-1}) \leq \Delta(\C,\eps).
\end{equation}
\end{thm}

\begin{remark} A more precise upper bound follows easily from the
proof, namely
$$
\Delta(\C,\eps) \le \left( \sum_{i\in \I}  (1+ \textup{deg}\
\phi_i) \right)\ \log(\eps^{-1}) + o(\log(\eps^{-1})).
$$ 
A more precise lower bound of the same kind can also be obtained for a
large class of Cantor sets which are generated by a set of full measure of
some random polynomial IFS's (see below for the definition of a random IFS). 
\end{remark}

\begin{remark}
The classical examples of Cantor sets are the middle $\frac 1
\rho$-th Cantor sets in the unit interval ($\rho=\frac 1 3$ gives
the usual middle third Cantor set). They can be thought of as
generated by IFS with affine contractions $\phi_0$ and $\phi_1$,
and contraction rate $\rho$. For these Cantor sets, all 
we need for their construction is the knowledge of $\eps$ and of $\rho$.
We can choose $\eps$ to be a number of low complexity.
In the particular case where $\rho$ is rational, $\Delta(\C,\eps)$ grows slower
than any unbounded partial recursive function.
\footnote{We thank Arnaldo Mandel for this observation.}
\end{remark}

\begin{remark}
Notice that there are examples of Cantor sets with low
$\eps$-distortion complexity which contain numbers of high
complexity. This happens for example in the middle third Cantor
set.
\end{remark}

Our next result is about real analytic IFS's.

\begin{thm} \label{thm:analytic}
Let $\C$ be a Cantor set generated by an IFS with real analytic
functions. Then
$$
\Delta(\C,\eps)\preccurlyeq \big(\log(\eps^{-1})\big)^2.
$$
\end{thm}

\bigskip

We now turn to random affine IFS's, for which at each step of the
construction we consider a random choice for the contraction rate.
For the definition we follow \cite{bamon}.

Let us consider a family $(\lambda_k)_{k\in \N}$ of
independent identically distributed random variables with values
in the interval $(0,1)$. To each sequence
$\lambda=(\lambda_k)_{k\in\N}$ we associate a Cantor set $\C$ in
the following way. Let $C^0_\lambda := [0,1]$. We define
$$
J^1_1(\lambda) := \left[0,\frac{\lambda_1}{2}\right],\
J^1_2(\lambda):= \left[1- \frac{\lambda_1}{2},1\right],
\ \textup{and}\
C^1_\lambda := J^1_1(\lambda) \cup J^1_2(\lambda) .
$$
In words, $C^1_\lambda$ is obtained by removing the central interval
of length $(1-\lambda_1)$ from $C^0_\lambda$. At the $(k+1)$-st step,
we delete from each interval $J^k_i(\lambda)$, $i=1,\dots,2^k$,
the central interval of length $(1-\lambda_{k+1})$, obtaining
$2^{k+1}$ intervals $J^{k+1}_i(\lambda)$, $i=1,\dots, 2^{k+1}$,
such that
\begin{equation} \label{rand-cantor-length}
|J^{k+1}_i(\lambda)|= \frac{1}{2^{k+1}} \prod_{h=1}^{k+1} 
\lambda_h \qquad \forall  i=1,\dots,2^{k+1}.
\end{equation}
Then we define
$$
C^{k+1}_\lambda := \bigcup_{i=1}^{2^{k+1}} J^{k+1}_i(\lambda).
$$
We call \emph{random central Cantor set} the set
$$
\C_\lambda := \bigcap_{k=0}^\infty C^k_\lambda.
$$
We remark that by construction the boundary points $\partial
J^k_i(\lambda)$ of all the intervals $J^k_i(\lambda)$ are
contained in $\C_\lambda$.

The next theorem states that random central Cantor sets need more information
than those generated by polynomial IFS's. 

\begin{thm} \label{thm-centcanrand}
Let $\C_\lambda$ be a random central Cantor set as described above.
Then, for any $\lambda \in (0,1)^\N$,
\begin{equation} \label{rand-1}
\Delta(\C_\lambda,\eps) \preccurlyeq
\left(\log(\eps^{-1})\right)^{2}.
\end{equation}
Moreover, let us assume that the common distribution of the i.i.d. random
variables $(\lambda_k)$ is absolutely continuous, with a density
$f(x)$ bounded above and below away from zero. Then, for any
$\delta >0$, we have
\begin{equation} \label{rand-2}
\left( \log(\eps^{-1})\right)^{2-\delta} \preccurlyeq
\Delta(\C_\lambda,\eps)
\end{equation}
for almost every $\lambda \in (0,1)^\N$.
\end{thm}

We now consider Cantor sets with a differentiable structure.
Following \cite{sullivan}, \cite{feliks} and \cite{bedford}, this
corresponds to Cantor sets generated by $C^k$ IFS's and we call
them $C^k$ Cantor sets. We shall recall their contruction in
Subsection \ref{sec:ck-central}.

\begin{thm}\label{thm:ckupper}
Let $k\ge 1$. For any $\delta>0$, for any $C^k$ Cantor set $\C$ with
box counting dimension $D$, we have
\begin{equation} \label{ck-1}
\Delta(\C,\eps)\preccurlyeq \eps^{-\frac D k -\delta}.
\end{equation}
Moreover, for any $\delta>0$, there exist (many) $C^{k}$ central
Cantor sets $\C$ with box counting dimension, at most $D+\delta$, such that
\begin{equation} \label{ck-2}
\eps^{-\frac D k+\delta} \preccurlyeq \Delta(\C,\eps).
\end{equation}
\end{thm}

We emphasise that in this case the asymptotic behaviour of the $\eps$-distortion complexity,
when $\eps$ tends to zero,  depends in general on the regularity $k$ of the set and, contrarily to the
previous cases, it also depends on its box counting dimension $D$.

\begin{remark}
We notice that our proofs also provide estimates for the Kolmogorov's $\eps$-entropy of some
families of Cantor sets \cite{titi} in the Hausdorff distance.
\end{remark}

\section{Proofs}\label{allproofs}

The following two simple lemmas will be used repeatedly in the proofs
hereafter. We leave their elementary proof to the reader.

\begin{lemma} \label{lemma:1}
Let $F$ and $F'$ be closed subsets of $A$. Let $I=[a,b]$ and
$I'=[a',b']$ be closed sub-intervals of $A$. Let $H=[c,d]$ and
$H'=[c',d']$ be closed subsets of $\stackrel{\scriptscriptstyle{\circ}}{I}$ and
$\stackrel{\scriptscriptstyle{\circ}}{I'}$ respectively. Assume that $\partial H
\subset F$, $\partial H' \subset F'$, $F\cap \stackrel{\scriptscriptstyle{\circ}}{H} =
\emptyset$ and $F'\cap \stackrel{\scriptscriptstyle{\circ}}{H'} = \emptyset$. Moreover
assume that there exists $\eps >0$ such that $|a-a'|\le \eps$,
$|b-b'|\le \eps$, $|c-d|>2\eps$, $|c'-d'|>2\eps$ and $\max
\set{|c-c'|,|d-d'|}>\eps$, then $\disth(F,F')>\eps$.
\end{lemma}

In the sequel, this lemma will be used to show that two Cantor sets ($F$ and $F'$)
are a Hausdorff distance larger than $\eps$, $H$ and $H'$ playing the role of
holes in the Cantor sets.

\begin{lemma} \label{lemma:2}
Let $(\Omega,\cal{A},\proba)$ be a probability space. Let $\C$
be a measurable map from $(\Omega,\cal{A})$ to the set of closed
subsets of $A$ equipped with the Borel $\sigma$-algebra induced by
the Hausdorff metric. Let $(a_k)_k$ be a positive, increasing, diverging
sequence. Assume that for any integer $k$ there exists a sequence
$(V_{k,j})_{1\le j\le 2^{a_k}}$ of measurable subsets of $\Omega$
such that
$$
\set{\omega : \Delta(\C(\omega),2^{-k}) < a_k} \subset
\bigcup_{j=1}^{2^{a_k}}\ V_{k,j}
$$
and
$$
\sum_{k}\ \sum_{j=1}^{2^{a_k}}\ \proba(V_{k,j}) < \infty.
$$
Then, for $\proba$-almost every $\omega$,
$\Delta(\C(\omega),2^{-k}) \ge a_k$ for any $k$ large enough
(depending on $\omega$).
\end{lemma}

\subsection{Proof of Theorem \ref{thm:pol}} 

Let $N$ be the largest degree of the polynomial functions
$\set{\phi_i}$, then we can write
$$\phi_i(x) = \sum_{0\le \alpha \le N} \ c_{i,\alpha} x^\alpha \qquad
\forall \ i \in \I$$ with coefficients $c_{i,\alpha} \in \R$. We
now show how to construct a program $\Pro$ approximating $\C$
within Hausdorff distance $\eps$.

Let $\eps$ be fixed and $K$ a constant to be specified later on. Let
us define $\eps'=\frac{\eps}{K}$. We construct polynomials
$$
\tilde \phi_i(x) := \sum_{0\le \alpha \le N} \ \tilde c_{i,\alpha}
x^\alpha \qquad \forall \ i \in \I
$$
with coefficients satisfying
\begin{equation} \label{polyn-app-coeff}
|c_{i,\alpha} - \tilde c_{i,\alpha}| < \eps' \qquad \forall\ i\in
\I \quad \forall\ 0\le \alpha \le N
\end{equation}
such that the $\tilde \phi_i$'s are injective contractions on $A$
with uniform contraction rate $\tilde \rho \in (0,1)$. For any
$\omega_1^n \in \I^n$ we construct the composition $\tilde
\phi_{_{\omega_1^n}}$ as in (\ref{ifs-finite}).

We first show that for any bounded set $B\subset \R$ such that
$\phi_i(B)\subset B$ for all $i\in \I$ with the same contraction
rate $\rho$, and $\tilde \phi_i(B)\subset B$ for all $i\in \I$, we
have for all $n\in \N$
\begin{equation} \label{polyn-estim-finite}
\disth(\phi_{_{\omega_1^n}}(B), \tilde \phi_{_{\omega_1^n}}(B)) <
\eps' \ \frac{1-\rho^n}{1-\rho}\ \sup\limits_{x\in B}\ \left|
\sum_{0\le \alpha \le N}\ x^\alpha \right| \qquad \forall\
\omega_1^n \in \I^n.
\end{equation}
The proof is by induction. By (\ref{polyn-app-coeff}) and
definition of $\disth$, one immediately gets
$$
\disth(\phi_i(B),\tilde \phi_i(B)) < \eps' \ \sup\limits_{x\in B}\
\left| \sum_{0\le \alpha \le N}\ x^\alpha \right| \qquad \forall \
i\in \I.
$$
The inductive step is obtained by using the triangle inequality
for $\disth$. By the first step we have
$$
\disth\left(\phi_{\omega_n}(\tilde \phi_{_{\omega_1^{n-1}}}(B)), \tilde
\phi_{\omega_n} (\tilde \phi_{_{\omega_1^{n-1}}}(B))\right) < \eps' \
\sup\limits_{x\in B}\ \left| \sum_{0\le \alpha \le N}\ x^\alpha
\right|
$$
where we have used $\tilde \phi_{_{\omega_1^{n-1}}}(B) \subset B$.
Moreover, by using the contraction rate $\rho$, we get
$$
\disth\left(\phi_{\omega_n}(\phi_{_{\omega_1^{n-1}}}(B)),
\phi_{\omega_n} (\tilde \phi_{_{\omega_1^{n-1}}}(B))\right) < \rho\
\disth(\phi_{_{\omega_1^{n-1}}}(B),\tilde
\phi_{_{\omega_1^{n-1}}}(B))
$$
$$
< \eps' \ \rho\ \frac{1-\rho^{n-1}}{1-\rho}\ \sup\limits_{x\in B}\
\left| \sum_{0\le \alpha \le N}\ x^\alpha \right|
$$
where the last inequality is the $(n-1)$-th step of the induction.
Hence the triangle inequality implies that
$$
\disth(\phi_{_{\omega_1^n}}(B), \tilde \phi_{_{\omega_1^n}}(B))
<\eps' \ \left(1+\rho\ \frac{1-\rho^{n-1}}{1-\rho}\right)\
\sup\limits_{x\in B}\ \left| \sum_{0\le \alpha \le N}\ x^\alpha
\right|.
$$
This finishes the proof of (\ref{polyn-estim-finite}).

Let us choose $\bar n\in \N$ such that $\rho^{\bar n} < \eps'$ and
$\tilde \rho^{\bar n} < \eps'$. For this fixed $\bar n$, let
$V:=\set{0,1}=\partial[0,1]$ and define
$$
\Ca := \bigcup_{\omega_1^{\bar n} \in \I^{\bar n}}\ \tilde
\phi_{_{\omega_1^{\bar n}}} (V).
$$
We now prove that $\disth(\C,\Ca) < \eps$. Let us consider $x\in
\C$ and $y\in \Ca$. By (\ref{ifs-prop}) there exists $z\in \C$
such that $x= \phi_{_{\omega_1^{\bar n}(x)}}(z)$ for a given
sequence $\omega_1^{\bar n}(x) \in \I^n$. Hence
\begin{equation} \label{pointwise-estim}
\begin{array}{c}
d(x,y)= d(\phi_{_{\omega_1^{\bar n}(x)}}(z),y) \le \\[0.5cm]
\disth(\phi_{_{\omega_1^{\bar
n}(x)}}(z),\phi_{_{\omega_1^{\bar n}(x)}}(V)) +
\disth(\phi_{_{\omega_1^{\bar n}(x)}}(V), \tilde
\phi_{_{\omega_1^{\bar n}(x)}}(V)) + \disth(\tilde
\phi_{_{\omega_1^{\bar n}(x)}}(V), y).
\end{array}
\end{equation}
For the first term we use the contraction properties to get
$$
\disth(\phi_{_{\omega_1^{\bar n}(x)}}(z),\phi_{_{\omega_1^{\bar n}(x)}}(V)) 
< \rho^{\bar n} \hbox{ diam } (A) < \eps'.
$$
By (\ref{polyn-estim-finite}), for the second term we have
$$
\disth(\phi_{_{\omega_1^{\bar n}(x)}}(V), \tilde
\phi_{_{\omega_1^{\bar n}(x)}}(V)) < \frac{\eps'}{1-\rho}\ (N+1).
$$
If we take
$$
K= 1 + \frac{N+1}{1-\rho}
$$
then
$$
d(x,y) \le \eps' K + \disth(\tilde \phi_{_{\omega_1^{\bar n}(x)}}(V), y)=
\eps + \disth(\tilde \phi_{_{\omega_1^{\bar n}(x)}}(V), y).
$$
Choosing $y\in \tilde \phi_{_{\omega_1^{\bar n}(x)}}(V) \subset
\Ca$, we have
$$
\disth(\tilde \phi_{_{\omega_1^{\bar n}(x)}}(V), y)=0,
$$
hence $d(x,y)< \eps$. Therefore
$$
\sup\limits_{x\in \C} \ d(x,\Ca) < \eps.
$$
On the other hand for a given $\omega_1^{\bar n} \in \I^{\bar n}$
and $y\in \tilde \phi_{_{\omega_1^{\bar n}}}(V)$, take $x\in
\phi_{_{\omega_1^{\bar n}}}(A) \cap \C$, noticing that this set is
not empty. Then we deduce that
$$
\sup\limits_{y\in \Ca} \ d(y,\C) < \eps.
$$
Hence $\disth(\C,\Ca) < \eps$.

Let us define the program $\Pro$ that contains the numbers $\eps$,
$\rho$ and $K$, and such that it specifies the coefficients
$\set{\tilde c_{i,\alpha}}$, computes $\bar n$ and makes the
computation of the $\tilde \phi_i(V)$'s. The binary length
$\ell(\Pro)$ satisfies
$$
\ell(\Pro) \preccurlyeq \log \left((\eps')^{-1} \right) =
\log K + \log(\eps^{-1})\ .
$$
Indeed, $\eps$ is specified with $O(\log \eps^{-1})$ bits, and $\rho$
and $K$ do not depend on $\eps$ and can be approximated by
rational numbers. The coefficients $\set{\tilde c_{i,\alpha}}$ are
approximations of the $\set{c_{i,\alpha}}$ with precision $\eps'$,
hence we can choose them as rational numbers requiring only
$O(\log (\eps')^{-1})$ bits of information. Finally the
information for the computation of $\bar n$ and $\Ca$ needs only
$O(1)$ bits of information. Hence this proves (\ref{pol-1}).

We now prove (\ref{pol-2}). Let $\I=\set{0,1}$ and define
$\phi_1(x)=bx$, $\phi_2(x)=1-bx$ for $b\in (0,1/2)$. We denote by
$\C_b$ the Cantor set generated by the IFS $\set{\phi_1,\phi_2}$.
To be in the context of Lemma \ref{lemma:2}, we take $b$ at random
according to the uniform distribution on the interval $(1/4,1/3)$.
We restrict the possible values of $b$ to ensure that the middle
hole is large enough so that an obviously simplified version of
Lemma \ref{lemma:1} applies.

For a fixed $\delta \in (0,1)$, define $a_k:=(1-\delta) k$. For
any $k$, there are at most $2^{a_k}$ different binary programs
$(\Pro_j)_{1\le j \le 2^{a_k}}$ of length $a_k-1$, which generate
at most $2^{a_k}$ different sets $\Ca_j:=\Ca(\Pro_j)$. We define
$$
V_{k,j}:= \set{b : \disth(\C_b,\Ca_j) < 2^{-k}}.
$$
Then
$$
\set{b : \Delta(\C_b,2^{-k})< a_k} \subset
\bigcup_{j=1}^{2^{a_k}}\ V_{k,j}.
$$
We now estimate $\proba(V_{k,j})$. We denote by $\partial^+ \Ca_j$
the rightmost point of $\Ca_j \cap (0,1/2)$. For $k$ large enough
($k\ge 3$) and for a given $j$, if $b\in V_{k,j}$ then
$d(b,\partial^+ \Ca_j)< 2^{-k}$ and therefore $\proba(V_{k,j}) <
2^{1-k}$. This implies that
$$
\sum_k \ \sum_{j=1}^{2^{a_k}}\ \proba(V_{k,j}) < \sum_k \ 2^{a_k}\
2^{1-k} = \sum_k \ 2^{1-\delta k} < \infty.
$$
The result follows from Lemma \ref{lemma:2}.

\begin{remark}
Notice that this proof works also in arbitrary finite dimension.
\end{remark}

\subsection{Proof of Theorem \ref{thm:analytic}} 

We give the proof in the case that the $\set{\phi_i}$ are analytic
functions on an open ball $B(0,R)\subset \Com$ of radius $R>1$.
The general case follows by applying the same argument to
piecewise polynomial approximations of the $\set{\phi_i}$.

By hypothesis we can write for $z\in A$
$$
\phi_i(z) = \sum_{h=0}^\infty\ c_{i,h} \ z^h \qquad \forall \ i\in
\I
$$
with coefficients $c_{i,h} \in \Com$. By the analyticity of the
functions $\set{\phi_i}$ in $B(0,R)$ it follows that
$$
\limsup\limits_{h\to \infty} \ |c_{i,h}|\ R^h \le 1\ .
$$
Hence there exists a real constant $r$ such that
$$
|c_{i,h}|\ R^h \le r \qquad \forall \ h\ge 0.
$$
As above let us denote by $V=\set{0,1}=\partial[0,1]$. We now construct
approximations of the analytic functions $\set{\phi_i}$.

Let $\delta\in (0,1)$ be such that $\delta R>1$, and, for $\eps$
fixed, define
\begin{equation} \label{analytic-help}
\eps':= \frac{r(\delta R-1)}{R(1-\delta)}\ ((1-\delta)
\eps)^{\left( \frac{\log R}{\log \frac{1}{\delta}} \right)}.
\end{equation}
Let $N$ be the smallest integer satisfying
\begin{equation} \label{analytic-app-degree}
\sum_{h=N}^\infty \ \delta^h = \frac{\delta^N}{1-\delta} <
(1-\rho)\ \frac{\eps}{4r} \cdot
\end{equation}
Hence we construct the polynomials
$$
\tilde \phi_i(z) := \sum_{h=0}^{N-1} \ \tilde c_{i,h} \ z^h \qquad
\forall\ i\in \I
$$
with coefficients $\tilde c_{i,h} \in \Com$ such that
\begin{equation} \label{analytic-app-coeff}
|c_{i,h} - \tilde c_{i,h} | < \eps' \qquad \forall\ i \in \I \quad
\forall \ h=0,\dots,N-1
\end{equation}
and the $\tilde \phi_i$'s are contractions functions on $A$ with
$\tilde \phi_i (A) \subset A$ for all $i\in \I$, and uniform
contraction rate $\tilde \rho \in (0,1)$.

Let us choose ${\bar n}\in \N$ such that $\rho^{\bar n} <
\frac{\eps}{4R}$ and $\tilde \rho^{\bar n} < \frac{\eps}{4R}$.
For this fixed ${\bar n}$, we define
$$
\Ca := \bigcup_{\omega_1^{\bar n} \in \I^{\bar n}}\
\tilde \phi_{_{\omega_1^{\bar n}}}(V) .
$$
We now prove that $\disth(\C,\Ca) < \eps$.
The proof will follow again by using (\ref{pointwise-estim}) and the
analogue of (\ref{polyn-estim-finite}).

For any bounded set $B\subset B(0,\delta R)$ such that $\phi_i(B)
\subset B$ and $\tilde \phi_i(B) \subset B$, we have for all $n\in
\N$
\begin{equation} \label{analytic-estim-finite}
\disth(\phi_{_{\omega_1^n}}(B), \tilde \phi_{_{\omega_1^n}}(B)) <
\eps \ \frac{1-\rho^n}{2} \qquad \forall\ \omega_1^n \in \I^n \ .
\end{equation}
The proof of (\ref{analytic-estim-finite}) is by induction as the
proof of (\ref{polyn-estim-finite}). We only show the first step.
By using (\ref{analytic-app-coeff}) and
(\ref{analytic-app-degree}), we obtain for all $z\in B(0,\delta
R)$
$$
|\phi_i(z)-\tilde \phi_i(z)| = \left| \sum_{h=0}^{N-1}\
(c_{i,h}-\tilde c_{i,h}) \ z^h + \sum_{h=N}^\infty \ c_{i,h} \ z^h
\right| 
$$
$$\le \eps' \sum_{h=0}^{N-1}\ \delta^h R^h + r \sum_{h=N}^\infty \
\delta^h < \eps' \frac{\delta^N R^N}{R\delta -1} + (1-\rho)\
\frac{\eps}{4r} \ r
$$
$$
< (1-\rho)\ \eps \left(\eps' \frac{R(1-\delta)
((1-\delta)\eps)^{\big(-\frac{\log R}{\log \delta^{-1}}\big)} }{4r(\delta
R -1)} + \frac{r}{4r}\right)  = (1-\rho)\ \frac{\eps}{2}\ ,
$$
where in the last equality we have used (\ref{analytic-help}).
Hence
$$
\disth(\phi_i(B),\tilde \phi_i(B)) < \eps \ \frac{(1-\rho)}{2}
\qquad \forall \ i\in \I.
$$
The inductive step is obtained as in the proof of
(\ref{polyn-estim-finite}) by the triangle inequality for
$\disth$.

We now write (\ref{pointwise-estim}) and, by repeating the
argument of the proof of Theorem \ref{thm:pol} and by
(\ref{analytic-estim-finite}), we obtain $\disth(\C,\Ca) < \eps$.

Let us define the program $\Pro$ which contains the numbers $\eps$,
$R$, $r$, $\delta$, $N$ and ${\bar n}$, and such that it specifies
the coefficients $\set{\tilde c_{i,h}}$ for $h=0,\dots,N-1$, and
makes the computation of the $\tilde \phi_i(V)$'s. The binary length
$\ell(\Pro)$ satisfies
$$
\ell(\Pro) \preccurlyeq N \log \left((\eps')^{-1} \right)
\preccurlyeq \left( \log (\eps^{-1}) \right)^2\ . 
$$
Indeed, $\eps$ is specified with $O(\log \eps^{-1})$ bits, and $\bar n \asymp \log
\eps^{-1}$. Hence it is specified by $o(\log \eps^{-1})$ bits and
$N\asymp \log \eps^{-1}$, hence it is specified by $o(\log
\eps^{-1})$ bits. The coefficients $\set{\tilde c_{i,h}}$ are
approximations of the $\set{c_{i,h}}$ with precision $\eps'$,
hence we can choose them as rational numbers requiring only
$O(\log (\eps')^{-1})$ bits of information, but there are $N$ of
them for each function $\phi_i$ (see (\ref{analytic-help}) and
(\ref{analytic-app-degree})), hence we need $O(N \log
(\eps')^{-1})$ bits of information, that is $O(\left( \log
(\eps^{-1}) \right)^2)$. Finally $R$, $r$ and $\delta$ do not
depend on $\eps$, hence the information for them and for the
computation of $\Ca$ need only $O(1)$ bits of information. Hence
$$\Delta(\C,\eps) \preccurlyeq \left( \log (\eps^{-1}) \right)^2$$
and the theorem is proved.

\subsection{Proof of Theorem \ref{thm-centcanrand}}

We first prove (\ref{rand-1}) by constructing an approximation
$\Ca$ of the set $\C$.

Let us consider a fixed sequence $\lambda \in (0,1)^\N$ and the
Cantor set $\C_\lambda$. For a fixed $\eps$, let $\bar n$ be given by
\begin{equation} \label{rand-cantor-app-degree}
\bar n := \min \set{n\in \N:\frac{1}{2^n}< \frac \eps 2}.
\end{equation}
Next, let us consider a sequence $\tilde \lambda=(\tilde \lambda_k)_{k\in \N}\in(0,1)^{\N}$
such that $|\lambda_1 - \tilde \lambda_1| < \eps$ and
\begin{equation} \label{rand-cantor-coeff}
|\lambda_k - \tilde \lambda_k| < 2^{k-2}\ \eps \qquad \forall\
k=2,\dots,\bar n\ .
\end{equation}
Then we define the approximation $\Ca$ of $\C_\lambda$ to be the finite
set
$$
\Ca := \partial C^{\bar n}_{\tilde \lambda} =
\bigcup_{i=1}^{2^{\bar n}}\ \partial J^{\bar n}_i(\tilde \lambda)
$$
where the sets $C^{\bar n}_{\tilde \lambda}$ and $J^{\bar
n}_i(\tilde \lambda)$ are constructed as specified in Section
\ref{definition}. We now prove that $\disth(\C_\lambda,\Ca)
<\eps$. To this aim we show that
\begin{equation} \label{rand-cantor-app-finite}
\disth(\partial J^k_1(\lambda),\partial J^k_1(\tilde \lambda)) <
\frac \eps 2
\end{equation}
for all $k=1,\dots,\bar n$. The same argument applies to all other
sets $J^k_i$. This is enough since it implies that, for any two
points $x\in \C_\lambda$ and $y\in \Ca$ in the analogous intervals
({\em i.e.},  $x\in J^{\bar n}_i(\lambda)$ and $y\in
\partial J^{\bar n}_i(\tilde \lambda)$ with the same index
$i=1,\dots,2^{\bar n}$),
\begin{eqnarray*}
d(x,y)  & \le & d(x,\partial J^{\bar n}_i(\lambda)) + \disth(\partial J^{\bar
n}_i(\lambda),\partial J^{\bar n}_i(\tilde \lambda)) +
d(\partial J^{\bar n}_i(\tilde \lambda),y)\\
& < & \frac{1}{2^{\bar n}} + \frac \eps 2 + 0 < \eps ,
\end{eqnarray*}
where for the first term we have used (\ref{rand-cantor-length})
and (\ref{rand-cantor-app-degree}), for the second term we have
used (\ref{rand-cantor-app-finite}), and for the third term we
have used the definition of $\Ca$. Hence $\disth(\C_\lambda,\Ca)<
\eps$. It remains to prove (\ref{rand-cantor-app-finite}). By
definition
$$
\partial J^k_1(\lambda) = \set{0, \frac{1}{2^k} 
\prod_{h=1}^k  \lambda_h} , \partial J^k_1(\tilde \lambda) =
\set{0, \frac{1}{2^k}  \prod_{h=1}^k  \tilde \lambda_h}.
$$
Hence it is enough to show that
\begin{equation} \label{rand-cantor-app-finite-help}
d\left(\frac{1}{2^k}  \prod_{h=1}^k  \lambda_h, \frac{1}{2^k} 
\prod_{h=1}^k  \tilde \lambda_h\right) < \frac \eps 2
\end{equation}
for all $k=1,\dots,\bar n$. By definition of $\tilde \lambda$, it
holds
$$
d\left( \frac 1 2 \lambda_1, \frac 1 2 \tilde \lambda_1 \right) <
\frac \eps 2.
$$ 
Assuming that (\ref{rand-cantor-app-finite-help}) holds for $k-1$, one gets 
$$
d\left(\frac{1}{2^k}  \prod_{h=1}^k  \lambda_h, \frac{1}{2^k} 
\prod_{h=1}^k  \tilde \lambda_h\right) 
$$
$$
\le \left|\frac{1}{2^{k-1}}  \prod_{h=1}^{k-1}  \lambda_h\right|\
\left|\frac{\lambda_k - \tilde \lambda_k}{2} \right| +
d\left(\frac{1}{2^{k-1}}  \prod_{h=1}^{k-1}  \lambda_h,
\frac{1}{2^{k-1}}  \prod_{h=1}^{k-1}  \tilde \lambda_h\right) \
\frac{\tilde \lambda_k}{2}
$$
$$
< \frac{1}{2^{k-1}} \ \frac{2^{k-2} \eps}{2} + \frac{\eps}{4} < \frac \eps 2,
$$
where we have used (\ref{rand-cantor-coeff}) for the first term,
the inductive hypothesis for the second term, and the fact that
$\lambda, \tilde \lambda \in (0,1)^\N$.

To finish the proof of (\ref{rand-1}) we have to estimate the
length of a binary program $\Pro$ producing the finite set $\Ca$.
The program $\Pro$ must contain the number $\eps$, the information
to compute $\bar n$, the contraction factors $(\tilde \lambda_k)$
for $k=1,\dots,\bar n$ and the instruction to compute
$\Ca$. The number of instructions for all the computations are $O(1)$ with
respect to $\eps$. The number $\eps$ is specified by $O(\log
\eps^{-1})$ bits of information and $\bar n \approx \log \eps^{-1}$
and each coefficient $\tilde \lambda_k$ needs $O(\log \eps^{-1})$
bits of information. Since there are $\bar n$ coefficients to be
specified, we find
$$
\ell(\Pro) \preccurlyeq \left( \log (\eps^{-1}) \right)^2
$$
hence (\ref{rand-1}) follows.

We now prove (\ref{rand-2}). First of all we identify a full
measure set of \virg{good} $\lambda$. By the hypothesis on the
density $f(x)$ of the common distribution of the random variables
$(\lambda_k)$, the following quantity is finite
$$
\gamma:= \int_0^1 \ \log (x)\ f(x) \ dx < 0.
$$
Notice that $\frac{e^{\gamma}}{2}$ can be interpreted as the
typical contraction rate, since products of many i.i.d. random variables $\lambda_k$
will be involved.

Given any $\eta >0$ we define, for all $n\in \N$,
\begin{equation} \label{rand-cantor-null-set}
\Lambda_n:= \set{\lambda \in (0,1)^\N  : \prod_{h=1}^k 
\lambda_h > e^{k(\gamma-\eta)} \ , \forall\
k=\lfloor \sqrt{n} \rfloor,\dots, n}\ .
\end{equation}
We remark that, since $\lambda_k \in (0,1)$ for all $k \ge 1$, if
$\lambda \in \Lambda_n$ then
\begin{equation} \label{rand-cantor-bel-help}
\prod_{h=1}^k  \lambda_h > 
e^{\sqrt{n}(\gamma-\eta)}\qquad \forall\ k=1,\dots,\lfloor \sqrt{n}
\rfloor -1.
\end{equation}
 
\begin{lemma} \label{lem:centcanrand}
Let us denote by $\Lambda_n^c$ the complement of  $\Lambda_n$ in
$(0,1)^\N$, then
$$
\sum_n\ \proba(\Lambda_n^c) < \infty\ , 
$$
\emph{i.e.}, almost every $\lambda \in (0,1)^\N$ belongs to
$\Lambda_n^c$ only for finitely many $n\in \N$.
\end{lemma}

\noindent {\it Proof.} We use the large deviation principle for
independent and identically distributed random variables (see,
{\emph{e.g.}, \cite{DZ}). It implies that for any fixed $\eta >0$
there exists a positive 
constant $C$ such that
$$
\lim\limits_{k\to \infty}\ \frac{1}{k} \ \log \proba
\set{\frac{\log(\lambda_1\times\cdots\times\lambda_k)}{k} - \gamma <
  -\eta}
= -C\ .
$$
Hence, for $n$ large enough and for all $0<C'<C$, we have the
estimate
$$
\proba(\Lambda_n^c) \le \sum_{k=\lfloor \sqrt{n} \rfloor}^n 
\proba \set{\prod_{h=1}^k  \lambda_h < e^{k(\gamma-\eta)}} \le e^{-C' \sqrt{n}}\
\sum_{k=\lfloor \sqrt{n} \rfloor}^n \ e^{-C'(k-\sqrt{n})}\ .
$$
Therefore
$$
\proba(\Lambda_n^c) \preccurlyeq e^{-C' \sqrt{n}}
$$
and the lemma follows by the Borel-Cantelli Lemma. \qed

For any $\eps$ we define
\begin{equation} \label{rand-cantor-bel-degree}
N(\eps) := \min \set{n\in \N : \left(2 e^{\eta-\gamma}\right)^n  \eps > \left( \log (\eps^{-1})\right)^{-2}} .
\end{equation}
We now consider a subset of $\Lambda_{N(\eps)}$. Let $q\in \N$ and
define
$$
\Psi_q:= \set{\lambda \in \Lambda_{N(2^{-q})} :
(1-\lambda_{k+1})\ \frac{\prod_{h=1}^k \ \lambda_h}{2^k} >
2^{2-q} , \forall\ k=0,\dots,N(2^{-q})}.
$$

\begin{lemma} \label{lem:rand-useful}
We have
$$
\sum_q\ \proba(\Psi_q^c) < \infty
$$
\emph{i.e.}, almost every $\lambda \in (0,1)^\N$ lies in
$\Psi_q^c$ only for finitely many $q\in \N$.
\end{lemma}

\noindent {\bf Proof.} By Lemma \ref{lem:centcanrand} and the
Borel-Cantelli Lemma, it is enough to prove that
$$
\sum_q \ \proba \left( \Psi_q^c \cap \Lambda_{N(2^{-q})} \right) <
\infty.
$$
First observe that if $0\le k \le \sqrt{N(2^{-q})}$ and
$$
(1-\lambda_{k+1}) > 2^{2-q}\ 2^k\
(e^{\eta-\gamma})^{\sqrt{N(2^{-q})}}
$$
then $\lambda$ satisfies
\begin{equation} \label{cond-rand}
(1-\lambda_{k+1})\ \frac{\prod_{h=1}^k \ \lambda_h}{2^k} >
2^{2-q}.
\end{equation}
Similarly, if $\sqrt{N(2^{-q})}\le k \le N(2^{-q})$ and
$$
(1-\lambda_{k+1}) > 2^{2-q}\ 2^k\ (e^{\eta-\gamma})^k\ , 
$$
then (\ref{cond-rand}) holds. Therefore
$$
\proba \set{\lambda \in \Lambda_{N(2^{-q})} \setminus \Psi_q}\ \le
\ 2^{3-q}\ 2^{\sqrt{N(2^{-q})}}\
e^{(\eta-\gamma)\sqrt{N(2^{-q})}}\ +
$$
$$
+\ 2^{2-q}\ 2^{N(2^{-q})+1}\ e^{(\eta-\gamma)(N(2^{-q})+1)} \le
\frac{O(1)}{q^2}
$$
which is summable over $q$. The lemma is proved. \qed

For a given $\eps$, let $\lambda \in \Psi_{\log
(\eps^{-1})}$ and define
$$
M_{\eps,N(\eps)}(\lambda):=\Psi_{\log (\eps^{-1})}\; \bigcap 
$$
\begin{equation} \label{rand-cantor-bel-closed}
\left\{\tilde \lambda :
\begin{array}{ll}
|\lambda_1 - \tilde \lambda_1| < 2 \eps & \\
|\lambda_k - \tilde \lambda_k| < 2^{k+1} \left( e^{\eta-\gamma}
\right)^{\sqrt{N}} \eps & \textup{if}\ k=2,\dots,
\lfloor \sqrt{N} \rfloor -1 \\[0.2cm]
|\lambda_k - \tilde \lambda_k| < \frac{2}{e^{\eta-\gamma}} \left(
2\ e^{\eta-\gamma} \right)^k  \eps & \textup{if}\ k=\lfloor
\sqrt{N} \rfloor, \dots, N
\end{array}
\right\}
\end{equation}
We now show that if $\tilde \lambda \not\in M(\lambda)$ then
$\disth(\C_\lambda, \C_{\tilde \lambda}) \ge \eps$. This follows
from Lemma \ref{lemma:1} and we now check the hypothesis to be
satisfied.

If $\tilde \lambda \in \Psi_{\log (\eps^{-1})} \setminus
M(\lambda)$ then one of the conditions in
(\ref{rand-cantor-bel-closed}) is violated. Following the notation
of Lemma \ref{lemma:1}, we start with $I=I'=[0,1]$. We take
$H=[\frac{\lambda_1}{2}, 1-\frac{\lambda_1}{2}]$ and
$H'=[\frac{\tilde \lambda_1}{2}, 1-\frac{\tilde \lambda_1}{2}]$. Then
$|H|>2\eps$ and $|H'|>2\eps$ since $\lambda$ and $\tilde \lambda$
are in $\Psi_{\log (\eps^{-1})}$. If $|\lambda_1 - \tilde
\lambda_1| \ge 2 \eps$, then
$\left|\frac{\lambda_1}{2}-\frac{\tilde \lambda_1}{2}\right| \ge \eps$,
and Lemma \ref{lemma:1} applies with $F=\C_\lambda$ and $F'=\C_{\tilde
\lambda}$, implying $\disth(\C_\lambda,\C_{\tilde \lambda})>\eps$.

Assume that, for some $k=2,\dots, \lfloor \sqrt{N} \rfloor -1$, all
conditions in (\ref{rand-cantor-bel-closed}) are satisfied up to
$k-1$ and condition $k$ is violated.
Either there is an $\ell<k$ such that 
$$
\left| \frac{1}{2^\ell}  \prod_{h=1}^\ell  \lambda_h -
\frac{1}{2^\ell}  \prod_{h=1}^\ell  \tilde \lambda_h \right|>\eps,
$$
in which case we define $\hat k$ to be the smallest such $\ell$.
Or, if no such $\ell$ exists,  
$$
\left| \frac{1}{2^k}  \prod_{h=1}^k  \lambda_h -
\frac{1}{2^k}  \prod_{h=1}^k  \tilde \lambda_h \right| \ge
\frac{\left| \lambda_k - \tilde \lambda_k \right|}{2^k}
\prod_{h=1}^{k-1}  \lambda_h  - \frac{\tilde \lambda_k}{2^k}
\left| \prod_{h=1}^{k-1}  \lambda_h - \prod_{h=1}^{k-1}  \tilde
\lambda_h \right| \ge
$$
$$
\ge 2\ \eps \left( e^{\eta-\gamma} \right)^{\sqrt{N}}
\prod_{h=1}^{k-1}  \lambda_h - \eps > \eps
$$
where we have used (\ref{rand-cantor-bel-help}) and the fact  that the leftmost
positive points up to the $(k-1)$-th step of the construction are $\eps$-close, and we set $\hat k=k$.

We will apply Lemma \ref{lemma:1} with $I=J^{{\hat k}-1}_1(\lambda)$ and
$I'=J^{{\hat k}-1}_1(\tilde \lambda)$. We take
$$
H=\left[\frac{1}{2^{\hat k}}  \prod_{h=1}^{\hat k}  \lambda_h, \left( 1
-\frac{\lambda_{\hat k}}{2}\right) \frac{1}{2^{{\hat k}-1}}  \prod_{h=1}^{{\hat k}-1}
 \lambda_h \right]
$$
and
$$
H'=\left[\frac{1}{2^{\hat k}}  \prod_{h=1}^{\hat k}  \tilde\lambda_h, \left( 1
-\frac{\tilde\lambda_{\hat k}}{2}\right) \frac{1}{2^{{\hat k}-1}}  \prod_{h=1}^{{\hat k}-1}
 \tilde\lambda_h \right]\ .
$$
We have $|H|>2\eps$ and $|H'|>2\eps$
since $\lambda$ and $\tilde \lambda$ are in $\Psi_{\log
(\eps^{-1})}$, and we can apply Lemma \ref{lemma:1}
which gives $\disth(\C_\lambda,\C_{\tilde \lambda})>\eps$.

The same argument applies if the $k$-th condition
with $k=\sqrt{N},\dots,N$ is violated, and all conditions up to $k-1$ are
satisfied.
If the leftmost positive points up to the $(k-1)$-th step of the construction are $\eps$-close, we get 
$$
\left| \frac{1}{2^k}  \prod_{h=1}^k  \lambda_h -
\frac{1}{2^k}  \prod_{h=1}^k \tilde \lambda_h \right| \ge 2
\eps \left( e^{\eta-\gamma} \right)^{k-1} \prod_{h=1}^{k-1} 
\lambda_h - \eps > \eps
$$
by definition (\ref{rand-cantor-null-set}) of $\Lambda_{N(\eps)}$ 
and using, as above, that all previous leftmost positive points are
$\eps$-close to each other. Again this implies that
$\disth(\C_\lambda,\C_{\tilde \lambda})>\eps$.

Let us now estimate the measure of the set $M_{\eps,N(\eps)}(\lambda)$.
By an easy computation based on the independence of the random variables
$(\tilde{\lambda}_k)$, we obtain that for all $\lambda \in
\Psi_{\log(\eps^{-1})}$
$$
\proba(M_{\eps,N(\eps)}(\lambda)) 
$$
$$
\leq \big(\max\limits_{x\in [0,1]} f(x)\big)^{N(\eps)}
(2\eps)^{N(\eps)} (e^{\eta-\gamma})^{-\sqrt{{N(\eps)}}} \
2^{\sum_{k=2}^{N(\eps)} \ k} \
(e^{\eta-\gamma})^{\sum_{k=\lfloor\sqrt{{N(\eps)}}\rfloor}^{N(\eps)} \ k} 
$$
\begin{equation}
\label{emme}
\leq \left( 2\ e^{\eta-\gamma} \right)^{-\frac{{N(\eps)}^2}{2} +
O({N(\eps)})} = e^{-O(1) (\log(\eps^{-1}))^2} ,
\end{equation}
where we have used the definition (\ref{rand-cantor-bel-degree}) of ${N(\eps)}$.
Note that this estimate is uniform in $\lambda\in \Psi_{\log(\eps^{-1})}$.

For a fixed $\delta \in (0,1)$, define $a_q:=q^{2-\delta}$. For
any $q$, there are at most $2^{a_q}$ different binary programs
$(\Pro_j)_{1\le j \le 2^{a_q}}$ of length $a_q-1$, which generate
at most $2^{a_q}$ different sets $\Ca_j:=\Ca(\Pro_j)$. We define
$$
V_{q,j}:= \set{\lambda : \disth(\C_\lambda,\Ca_j) < 2^{-q}}.
$$
Then
$$
\set{\lambda  : \Delta(\C_\lambda,2^{-q})< a_q} \subset
\bigcup_{j=1}^{2^{a_q}}\ V_{q,j}.
$$
We can write
$$
\proba\left( \bigcup_{j=1}^{2^{a_q}} \ V_{q,j} \right) \le
\proba(\Lambda^c_{N(2^{-q})}) + \proba(\Psi_q^c \cap
\Lambda_{N(2^{-q})}) + \sum_{j=1}^{2^{a_q}}\ \proba(V_{q,j} \cap
\Psi_q).
$$
Moreover if $V_{q,j} \cap \Psi_q\not= \emptyset$, there is a
$\lambda \in \Psi_{q}$ such that $V_{q,j} \cap \Psi_q \subset
M_{2^{-q},N(2^{-q})}(\lambda)$. By Lemmas \ref{lem:centcanrand} and
\ref{lem:rand-useful}, and by (\ref{emme}) it follows that
$$
\sum_q \sum_{j=1}^{2^{a_q}}\ \proba(V_{q,j}) < \infty.
$$
The result follows from Lemma \ref{lemma:2}.

\subsection{Proof of Theorem \ref{thm:ckupper}}
\label{sec:ck-central}

{\bf Preliminaries.} We first recall the definition of the
\emph{scaling function} $S_\C$ of a $C^k$ Cantor set $\C$
(\cite{sullivan},\cite{feliks}). In the sequel we fix
$\I=\set{0,1}$. For a word $\omega_1^n \in \I^n$ we let
$$
J_{\omega_1^n}:= \phi_{_{\omega_1^n}}([0,1])=\phi_{\omega_n} \circ
\phi_{\omega_{n-1}} \circ \dots \circ \phi_{\omega_1}([0,1])\ .
$$
Then by definition
$\phi_i(J_{\omega_1^n})= J_{\omega_1^n i}$ for any $i\in \I$, and
it holds
$$
J_{\omega_1^n} \subset J_{\omega_2^n} \subset \dots \subset
J_{\omega_{n-1}^n}\subset J_{\omega_n}.
$$
The scaling function describes the contraction rates in the
previous inclusions. For a word $\omega_1^n \in \I^n$ we define
$\tilde S_\C(\omega_1^n) \in (0,1)^2$. The two components of
$\tilde S_\C(\omega_1^n)$ are the rates of contractions
$$
(\tilde S_\C(\omega_1^n))_i =\
\frac{|J_{i\omega_1^n}|}{|J_{\omega_1^n}|}\qquad i=0,1
$$
where $|J|$ denotes the length of the interval $J$. The length of
the gap between the two intervals $J_{0\omega_1^n}$ and
$J_{1\omega_1^n}$ in $J_{\omega_1^n}$ can be reconstructed from
these data. The scaling function is defined to be the
function
$$
S_\C : \I^\infty \to (0,1)^2
$$
given by
$$
S_\C (\omega) := \lim_{n\to \infty}\ \tilde
S_\C(\omega_1^n).
$$
We refer to \cite{sullivan} for the proof of the existence of this limit.

By definition one has
$$
|J_{\omega_1^n}| = \prod_{j=1}^{n-1} \left( \tilde S_\C
(\omega_{j+1}^n) \right)_{\omega_j}\ .
$$
By using the scaling function we can introduce a distance
$d_S(\omega,\tilde \omega)$ on $\I^\infty$ in the following way.
For two sequences $\omega,\tilde \omega \in \I^\infty$, let
$\omega \cap \tilde \omega$ denote their longest common prefix,
and let $|\omega \cap \tilde \omega|$ denote its length. Then we
let
\begin{equation} \label{scaling-distance}
d_S(\omega,\tilde \omega) := \sup\limits_{\alpha \in I^\infty}\
\prod_{j=1}^{n=|\omega \cap \tilde \omega|} \left( S_\C
(\omega_{j+1}^n \ \alpha) \right)_{\omega_j}\ .
\end{equation}
Then there exists a constant $H>0$ such that for any
$\omega,\tilde \omega \in \I^\infty$ it holds
$$
\frac 1 H \le \frac{J_{\omega \cap \tilde
\omega}}{d_S(\omega,\tilde \omega)}\le H \ .
$$
Relations between the properties of the scaling function of a
$C^k$ central Cantor set and the differentiability of the IFS
generating this Cantor set have been studied in \cite{feliks} in
the case $k\ge 1$ (we refer the reader to Main Theorem
\cite{feliks}, page 406). The idea is the following. Let
$A(\omega_1^n)$ denote the set of the four boundary points of the
intervals $(J_{i\omega_1^n})_{i\in \I}$. A scaling function $S_\C$
generates a $C^k$ Cantor set $\C$ if and only if for any $n\in \N$
there are diffeomorphisms from $A(\omega_1^n)$ into $A(\tilde
\omega_1^n)$, for any $\omega_1^n \not= \tilde \omega_1^n \in
\I^n$, with derivatives bounded by a constant $C(\omega_1^n,
\tilde \omega_1^n)$ which satisfies
\begin{equation} \label{ck-condition}
C(\omega_1^n, \tilde \omega_1^n) = C\ d_S(\omega_1^n \alpha,
\tilde \omega_1^n \alpha)^{k-1}
\end{equation}
where $C$ does not depend on $n$ and from the definition of $d_S$
the right hand side is independent on $\alpha \in \I^\infty$.

\noindent \textbf{Proof of Theorem \ref{thm:ckupper}.} We first
prove (\ref{ck-1}). Let $\eps$ be fixed. We show how to
approximate the set $\C$ within Hausdorff distance $\eps$.
We will  give the proof for integer $k\geq 1$. The proof easily
extends to functions whose $k$-th derivative is H\"older.

We can write the Taylor expansions of the maps $\phi_i$ at a point $x_0$
$$
\phi_i (x) = \sum_{p=0}^{k-1}\ c_{i,p}(x_0)\ (x-x_0)^p + R_i(x,x_0)
\qquad \forall\ x\in [0,1].
$$
Moreover there exists a constant $K>0$ such that $|R_i(x,x_0)|\le
K |x-x_0|^k$ for $x\in [0,1]$, for all $i\in \I$ and $x_0 \in
[0,1]$.

Let $\eps'= \eps \ \frac{1-\rho}{M}$ for a constant $M$ to be
specified later on. We now construct a sequence of polynomials which
approximate the maps $(\phi_i)_{i\in \I}$. If $D$ is the box
counting dimension of the Cantor set $\C$, we need for any
$\delta>0$ at most $N =O((\eps')^{-\frac D k-\delta})$ intervals
$(I_s)_{s=1,\dots,N}$ of size $(\eps')^{\frac 1 k}$ to cover $\C$.
Hence we can consider the maps $(\phi_i)_{i\in \I}$ restricted to
the sets $(I_s)$. If $y_s$ denotes the middle point of the
interval $I_s$, let $\tilde y_s$ be the approximation of the point
$y_s$ within a distance $\eps'$. Then we define
$$\tilde \phi_i^s (x) = \sum_{p=0}^k\ \tilde
c_{i,p}(y_s)\ (x-\tilde y_s)^p \qquad \forall\ x\in [0,1]$$ such
that
\begin{equation} \label{ck-app-coeff}
|c_{i,p}(y_s) - \tilde c_{i,p}(y_s)| < \eps' \qquad \forall \ i\in
\I \ \ \forall\ p=0,\dots,k
\end{equation}
and they are contractions on $\R$ with the same uniform
contraction rate $\tilde{\rho} <1$.

To construct an approximation of $\C$, we work on the boundary
points of the intervals $J_{\omega_1^n}$ which all are in $\C$.
Let us denote
$J_{\omega_1^n}=[y^1_{\omega_1^n},y^2_{\omega_1^n}]$. Since for
any $n\in \N$ and any $\omega_1^n \in \I^n$ we have
$y^\eta_{\omega_1^n}\in \C$ for $\eta=1,2$, we can associate to a
given $y^\eta_{\omega_1^n}$ a sequence $\sigma_0^{n-1} \in
\set{1,\dots,N}^{n-1}$ which specifies to which intervals of the
cover $(I_s)_{s=1,\dots,N}$ the pre-images
$y^\eta_{\omega_1^{n-1}}, y^\eta_{\omega_1^{n-2}}, \dots,
y^\eta_{\omega_1}, y^\eta_{\sharp}$ of $y^\eta_{\omega_1^n}$
belong, where $y^\eta_{\sharp} \in \set{0,1}$.

We now establish the analogue of (\ref{polyn-estim-finite}) for
the boundary points. Let us define
$$
\tilde y^\eta_{\omega_1^n}:= \tilde \phi_{\omega_n}^{\sigma_{n-1}}
\circ \tilde \phi_{\omega_{n-1}}^{\sigma_{n-2}} \circ \dots \circ
\tilde \phi_{\omega_1}^{\sigma_0} (y^\eta_\sharp)
$$
then for all $n\in \N$ it holds
\begin{equation} \label{ck-estim-finite}
|y^\eta_{\omega_1^n} - \tilde y^\eta_{\omega_1^n}| < \left( k+1 +
K + \max\limits_{i\in \I,\ s=1,\dots,N} \sum_{p=0}^k\
|c_{i,p}(y_s)| \right) \frac{\eps'}{1-\rho} \qquad \forall\
\omega_1^n \in \I^n.
\end{equation}
The proof is by induction. The first step ($n=1$) follows by
definition of the approximating polynomials,  estimates
(\ref{ck-app-coeff}) and properties of the remainder
$R_i(x,x_0)$. This yields
\begin{eqnarray*}
|y^\eta_i-\tilde y^\eta_i| & =  & |\phi_i(y^\eta_\sharp)- \tilde \phi_i^{\sigma_0}(y^\eta_\sharp)| \\
& \le & \eps'\sum_{p=0}^k (\eps')^{\frac p k} + \sum_{p=0}^k\ |c_{i,p}(y_{\sigma_0})|
|(y^\eta_\sharp-y_{\sigma_0})^p-(y^\eta_\sharp-\tilde
y_{\sigma_0})^p| + K \eps' \\
& \le & \left( k+1 + K + \max\limits_{i\in \I,\ s=1,\dots,N} \sum_{p=0}^k\ |c_{i,p}(y_s)|
\right) \eps'  .
\end{eqnarray*}
The inductive step follows by using the triangle inequality
$$
|y^\eta_{\omega_1^n} - \tilde y^\eta_{\omega_1^n}|\le
|\phi_{\omega_n}(y^\eta_{\omega_1^{n-1}})-\tilde
\phi_{\omega_n}^{\sigma_{n-1}}(y^\eta_{\omega_1^{n-1}})| + |\tilde
\phi_{\omega_n}^{\sigma_{n-1}}(y^\eta_{\omega_1^{n-1}})-\tilde
\phi_{\omega_n}^{\sigma_{n-1}}(\tilde y^\eta_{\omega_1^{n-1}})|
$$
together with
$$
|\phi_{\omega_n}(y^\eta_{\omega_1^{n-1}})-\tilde
\phi_{\omega_n}^{\sigma_{n-1}}(y^\eta_{\omega_1^{n-1}})| < \eps'
\left( k+1 + K + \max\limits_{i\in \I,\ s=1,\dots,N} \sum_{p=0}^k\
|c_{i,p}(y_s)| \right)
$$
and
$$
|\tilde \phi_{\omega_n}^{\sigma_{n-1}}(y^\eta_{\omega_1^{n-1}})-\tilde
\phi_{\omega_n}^{\sigma_{n-1}}(\tilde y^\eta_{\omega_1^{n-1}})| <
\rho\ |y^\eta_{\omega_1^{n-1}} - \tilde y^\eta_{\omega_1^{n-1}}|, 
$$
where $\rho$ is the uniform contraction rate of the approximating
polynomials.

Let us choose $\bar n$ such that $\rho^{\bar n} < \frac{\eps}{2}$
for all $n\ge \bar n$. Then we define the set
$$\Ca := \bigcup_{\omega_1^{\bar n}\in \I^{\bar n}}\ \left(\tilde
y^1_{\omega_1^{\bar n}}\cup \tilde y^2_{\omega_1^{\bar n}}
\right)$$ and we claim that $\disth(\C,\Ca)< \eps$. Indeed, by
definition of $\bar n$, any point in the Cantor set $\C$ is at
most at distance $\frac \eps 2$ from a point in the boundary of
one of the sets $J_{\omega_1^{\bar n}}$. Moreover, by construction
of the points $\tilde y^\eta_{\omega_1^{\bar n}}$ we have
(\ref{ck-estim-finite}), hence the claim follows since
$$
\disth(\C,\Ca) \le \sup\limits_{x\in \C}
\inf\limits_{\omega_1^{\bar n} \in \I^{\bar n}} |x-\tilde
y^\eta_{\omega_1^{\bar n}}| 
$$
$$
\le \sup_{x\in \C}
\inf\limits_{\omega_1^{\bar n} \in \I^{\bar n}} \left( |x-
y^\eta_{\omega_1^{\bar n}}| + |y^\eta_{\omega_1^{\bar n}} - \tilde
y^\eta_{\omega_1^{\bar n}}| \right)< \eps
$$
provided that we choose
$$
M:= 2 \left(
k+1 + K + \max\limits_{i\in \I,\ s=1,\dots,N} \sum_{p=0}^k\
|c_{i,p}(y_s)| \right)\ .
$$

Let us define the program $\Pro$ that contains the numbers $\eps$,
$\rho$, $D$, $M$, $K$ and $k$, and such that it specifies all the
necessary coefficients $\tilde c_{i,p}$, makes the computation to
obtain $N$ and the approximated points $\tilde y_s$, and moreover
it makes the computations to obtain $\bar n$ and the points
$\tilde y^\eta_{\omega_1^{\bar n}}$. The binary length
$\ell(\Pro)$ satisfies
$$
\ell(\Pro) \preccurlyeq \ \eps^{-\frac D k-\delta}
$$
since $\eps$ is specified with $O(\log(\eps^{-1}))$ bits; $\rho$,
$D$, $M$, $K$ and $k$ do not depend on $\eps$ and can be
approximated by rational numbers. The coefficients $\tilde
c_{i,p}$ and the points $\tilde y_s$ are constructed as in the
previous proofs with precision $\eps'$, hence each of them needs
$O(\log(\eps^{-1}))$ bits of information and their number is $N(k+2) =
O(\eps^{-\frac D k-\delta})$. Finally all the computations to
obtain $\Ca$ need $O(1)$ bits of instructions. Hence (\ref{ck-1})
follows.

We now prove (\ref{ck-2}). We define a class of particular scaling
functions $S(\alpha)$ on $\I^\infty$ to construct differentiable
Cantor sets with the given distortion complexity.

Let us denote by $\I^*:= \cup_{n\in \N}\ \I^n$ the countable set
of finite strings $s$ written using the alphabet $\I$. Let
$(\lambda_s)_{s \in \I^*}$ be a family of independent
identically distributed random variables with values in the
interval $(0,1)$ and absolutely continuous distribution with
density $f(x)$ bounded above and below away from zero. Note that the
empty string $\sharp$ belongs to $\I^*$ and therefore there is an
associated random variable $\lambda_\sharp$.

Let $0<\zeta <1$, $0<\rho <1$ and $\rho<\theta<1$ be given
constants, with $\rho$ determining the contraction rate. $\zeta$ will be chosen
small enough later on. We will
only consider central Cantor sets, namely the two components of
the scaling functions will be equal. We define the scaling
function
$$
S_\lambda(\alpha) := \rho + \zeta \sum_{q=1}^\infty\ \theta^{q-1}\
\lambda_{\alpha_1^{q}} \qquad \forall\ \alpha \in \I^\infty
$$
which depends on the realisation of the family $(\lambda_s)$.
We remark that for any realisation it holds
\begin{equation} \label{ck-scaling-estimate}
\rho \le S_\lambda(\alpha) \le \rho + \frac{\zeta}{1-\theta}
\qquad \forall\ \alpha \in \I^\infty\ .
\end{equation}
Hence if $\zeta$ is small enough, the rate of contraction is
almost $\rho$. It is also useful to define the truncated scaling
function $\tilde S_\lambda$ by
\begin{equation} \label{truncated}
\tilde S_\lambda(\omega_1^m) := \rho + \zeta \sum_{q=1}^{m} \
\theta^{q-1} \ \lambda_{\omega_1^{q}}  \quad \forall\ m\ge 0.
\end{equation}
Using the relations
\begin{equation} \label{ck-cantor}
\frac{|J_{i\omega_1^m}|}{|J_{\omega_1^m}|} = \tilde
S_\lambda(\omega_1^m) \qquad \forall \ i\in \I
\end{equation}
we can construct a central Cantor set $\C_\lambda$ generated by
the scaling function $S_\lambda$. From (\ref{ck-scaling-estimate})
it follows that the Cantor set $\C_\lambda$ has box counting
dimension $D(\zeta)$ which satisfies
\begin{equation} \label{ck-box-count}
D(\zeta) = -\ \frac{\log 2}{\log \rho} + O(\zeta) \qquad \mbox{ as
}\ \zeta \to 0.
\end{equation}

We now consider the differentiability of the IFS generating
$\C_\lambda$. By (\ref{ck-scaling-estimate}) and the definition
(\ref{scaling-distance}) of $d_S$ it follows
$$
\rho^n \le d_S(\omega,\tilde \omega) \le \left(
\rho + \frac{\zeta}{1-\theta} \right)^n \qquad n=|\omega \cap
\tilde \omega|
$$
for any $\omega,\tilde \omega \in \I^\infty$. Moreover for any
$\omega \not= \tilde \omega$ it holds
$$
|S_\lambda(\omega) - S_\lambda(\tilde \omega)| \le 2 \zeta \
\frac{\theta^{n+1}}{1-\theta} \qquad n=|\omega \cap \tilde
\omega|.
$$
Hence for $m>n$ we have diffeomorphisms from $A(\omega_1^m)$ into
$A(\tilde \omega_1^m)$ with derivatives bounded by a constant
$C(\omega_1^m, \tilde \omega_1^m) =O(\theta^{n+1})$. These facts
together with relation (\ref{ck-condition}) imply that the Cantor
set $\C_\lambda$ is of class $C^k$ with
\begin{equation} \label{ck-diff}
k = 1+ \frac{\log \theta}{\log \rho} \ge 1.
\end{equation}

Let $0<\eps<1$ be fixed and small enough depending on the
constants $\rho,\theta,\zeta$. Let $\lambda:=(\lambda_s)$ and
$\lambda':=(\lambda'_s)$ denote two different realisations of
the family of random variables. We give a condition on $\lambda$
and $\lambda'$ to have $\disth(\C_\lambda,\C_{\lambda'})> \eps$.
We denote
\begin{equation} \label{p-bar}
\bar p=\left[\frac{\log (C \eps\ \zeta^{-1})}{\log(\rho \theta)}\right]
\end{equation}
where $C$ is a positive constant (independent of $\eps$) to be specified later on.

For any $\sigma \in \I^*$, we denote by $J_{\sigma}$ and
$J'_{\sigma}$ the intervals associated to $\sigma$ in
$\C_{\lambda}$ and $\C_{\lambda'}$ respectively.

\begin{lemma} \label{lem:ck-dist}
Assume there is $0\le p\le \bar p$ satisfying
\begin{equation} \label{ipot}
\max\limits_{\omega_1^p\in \I^p}\
\disth(J_{\omega_1^p},J'_{\omega_1^p})
> \eps.
\end{equation}
Then $\disth(\C_\lambda,\C_{\lambda'})> \eps$.
\end{lemma}

\noindent \emph{Proof.} Denote by $p$ the smallest integer for
which the above inequality holds and by $\omega_1^p \in \I^p$ the
string realising the maximum. We apply Lemma \ref{lemma:1} with
$I=J_{\omega_1^{p-1}}$ and $I'=J'_{\omega_1^{p-1}}$. The
hypotheses on $I$ and $I'$ follow by the fact that (\ref{ipot}) is
violated up to $p-1$. The gaps $H$ and $H'$ have size at least
$$
\rho^{p-1} >(\rho \theta)^{\bar p} = C \eps\ \zeta^{-1}>2\eps
$$
for $\eps$ small enough if $C\zeta^{-1}>2$. 

Finally since $\disth(J_{\omega_1^p},J'_{\omega_1^p}) > \eps$ we have all the
hypotheses of Lemma \ref{lemma:1}. Hence the lemma follows. \qed

\begin{lemma} 
Assume that there exists $0\le p \le \bar p$ and a sequence
$\omega_1^p\in \I^p$ such that
$$
|\lambda_{\omega_1^p} - \lambda'_{\omega_1^p}| > (\rho
\theta)^{-p} \ \frac{(4+2\rho^2) \eps}{\zeta}\cdot
$$
Then $\disth(\C_{\lambda},\C_{\lambda'}) >\eps$.
\end{lemma}

\noindent \emph{Proof.} Denote by $p$ the smallest integer for
which the above inequality holds. It is enough to assume that for
any $q<p$ condition (\ref{ipot}) is not verified, otherwise the proof follows
immediately from Lemma \ref{lem:ck-dist}.

Let $x$ (respectively $x'$) be the point in the boundary of
$J_{\omega_1^p}$ (respectively $J'_{\omega_1^p}$) which is not in
the boundary of $J_{\omega_2^p}$ (respectively $J'_{\omega_2^p}$).
Let $y$ (respectively $y'$) be the other boundary point of
$J_{\omega_2^p}$ (respectively $J'_{\omega_2^p}$). Since by the
recursive assumption $|y-y'|\le \eps$, we have
$$
|x-x'|\ge |x-y-(x'-y')|-\eps.
$$
Now since $|x-y|=|J_{{\omega_1^p}}|$ and $(x-y)$ and $(x'-y')$
have the same sign, we get
$$
|x-x'|\ge \big| |J_{{\omega_1^p}}|-|J'_{{\omega_1^p}}|\big|-\eps
$$
This can also be written
$$
|x-x'|\ge \big|\tilde S_\lambda({\omega_2^p}) |J_{\omega_2^p}| -
\tilde S_{\lambda'}({\omega_2^p}) |J'_{\omega_2^p}| \big|- \eps
\ge
$$
$$
\ge \big|\tilde S_\lambda({\omega_2^p}) - \tilde
S_{\lambda'}({\omega_2^p}) \big|\
\frac{|J_{\omega_2^p}|+|J'_{\omega_2^p}|}{2} -\frac{\tilde
S_\lambda({\omega_2^p})+ \tilde S_{\lambda'}({\omega_2^p})}{2}\
\big| |J_{\omega_2^p}|-|J'_{\omega_2^p}|\big|-\eps \ge
$$
$$
\ge \big|\tilde S_\lambda({\omega_2^p}) - \tilde
S_{\lambda'}({\omega_2^p}) \big|\
\frac{|J_{\omega_2^p}|+|J'_{\omega_2^p}|}{2}\ -3 \eps \ge \rho^p
\big|\tilde S_\lambda({\omega_2^p}) - \tilde
S_{\lambda'}({\omega_2^p}) \big|\ - 3\eps
$$
since $\disth \left( J_{{\omega_2^p}}, J'_{{\omega_2^p}}\right)\le
\eps$ and $S_\lambda\ge \rho$. 
From \eqref{truncated} we get
$$
\tilde S_\lambda({\omega_2^p}) - \tilde S_{\lambda'}({\omega_2^p})
= \zeta \sum_{q=2}^{p} \theta^{q-2} \left(\lambda_{\omega_2^{q}}-\lambda'_{\omega_2^{q}} \right)
$$
and from (\ref{ck-cantor}) we have
$$
\zeta \sum_{q=2}^{p-1} \theta^{q-2} \left(\lambda_{\omega_2^{q}}-\lambda'_{\omega_2^{q}} \right) =
\frac{ (|J_{\omega_1^{p-1}}|-|J'_{{\omega_1^{p-1}}}|)
(|J_{{\omega_2^{p-1}}}|+|J'_{{\omega_2^{p-1}}}|)}{2\
|J_{{\omega_2^{p-1}}}||J'_{{\omega_2^{p-1}}}|} 
$$
$$
+\frac{ (|J_{{\omega_1^{p-1}}}|+|J'_{{\omega_1^{p-1}}}|)
(|J'_{{\omega_2^{p-1}}}|-|J_{{\omega_2^{p-1}}}|)}{2\
|J_{{\omega_2^{p-1}}}||J'_{{\omega_2^{p-1}}}|}\cdot
$$
Since $|J_{{\omega_1^{p-1}}}| \le |J_{{\omega_2^{p-1}}}| \le (\rho
+ \frac{\zeta}{1-\theta})^{p-2}$ we get (using again $\disth
\left( J_{{\omega_2^p}}, J'_{{\omega_2^p}}\right)\le \eps$,
and $S_\lambda \ge \rho$)
$$
\left|\ \zeta \sum_{q=2}^{p-1}\ \theta^{q-2} \left(
\lambda_{\omega_2^{q}}-\lambda'_{\omega_2^{q}} \right) \right|
\le 2\ \rho^{-p+2}\ \eps\ .
$$
Hence
$$
|\tilde S_\lambda({\omega_2^p}) - \tilde
S_{\lambda'}({\omega_2^p})| \ge \zeta \theta^{p-2} |\lambda_{{\omega_1^p}}-\lambda'_{{\omega_1^p}}|-2\ \rho^{-p+2}\
\eps .
$$
We conclude that
$$
|x-x'| \ge \zeta(\rho \theta)^p \
|\lambda_{{\omega_1^p}}-\lambda'_{{\omega_1^p}}|-
(3+2\rho^2)\eps.
$$
Therefore the lemma follows by applying Lemma \ref{lemma:1}. \qed

We now want to estimate for a given realisation
$\lambda$ of the family $(\lambda_s)$ the probability of the
event
$$
\mathscr{E}_{\bar p}(\lambda) = \sup_{0\le q\le {\bar p}} (\rho
\theta)^q\ \sup_{\sigma \in \I^{q}}
|\lambda_\sigma-\lambda'_\sigma| \le \frac{(4+2\rho^2)\
\eps}{\zeta}\cdot
$$
By independence of the family $(\lambda_s)$ we get 
$$
\proba \left( \mathscr{E}_{\bar p} (\lambda) \right) \le
\prod_{q=0}^{\bar p} \left(C\ (\rho \theta)^{-q}\ \eps \
\zeta^{-1} \right)^{2^{q}}
$$
where 
$$
C=(4+2\rho^2) \sup_{x\in[0,1]} f(x)
$$
where $f$ is the density of the random variables $(\lambda_s)$.

These relations imply that
\begin{eqnarray*}
\log \proba \left( \mathscr{E}_{\bar p} (\lambda) \right) & \le &
\sum_{q=0}^{\bar p}\ 2^{q} \left( \log (C \eps\ \zeta^{-1})
-q\log(\rho\theta) \right) \\
& \leq &  2^{\bar p +1}\log(\rho\theta) - (3+{\bar p}) \log(\rho\theta).
\end{eqnarray*}
For a fixed $\delta \in (0,1)$, define
$a_\ell:=2^{\ell(-\delta+\frac D k)}$ and choose $\eps =
2^{-\ell}$. For any $\ell$, there are at most $2^{a_\ell}$
different binary programs $(\Pro_j)_{1\le j \le 2^{a_\ell}}$ of
length $a_\ell-1$, which generate at most $2^{a_\ell}$ different
sets $\Ca_j:=\Ca(\Pro_j)$. We define
$$
V_{\ell,j}:= \set{\lambda : \disth(\C_\lambda,\Ca_j) < 2^{-\ell}}.
$$
Then
$$
\set{\lambda : \Delta(\C_\lambda,2^{-\ell})< a_\ell} \subset
\bigcup_{j=1}^{2^{a_\ell}}\ V_{\ell,j}.
$$
If $V_{\ell,j}$ is not empty, there exists $\lambda$ such that
$V_{\ell,j} \subset \mathscr{E}_{\bar p} (\lambda)$. Then
$$
\proba\left( \bigcup_{j=1}^{2^{a_\ell}} V_{\ell,j} \right) \le
2^{a_\ell} \proba\left( \mathscr{E}_{\bar p} (\lambda)\right) \le$$
$$
\le O(1) \exp \left( 2^{\ell(-\delta+\frac D k)} \log 2 +
2^{\bar p +1} \log(\rho\theta) - (3+{\bar p}) \log(\rho\theta)\right)\ .
$$
Using (\ref{ck-box-count}), (\ref{ck-diff}) and (\ref{p-bar}), it
follows that for any $\delta >0$ one can find $\zeta>0$ small
enough such that this is summable in $\ell$. Hence we can apply
Lemma \ref{lemma:2} to complete the proof. \qed

\end{document}